\documentclass[10pt]{article}
\usepackage[a4paper, margin=1in]{geometry}
\usepackage[utf8]{inputenc}
\usepackage[T1]{fontenc}
\usepackage{amsmath, amsthm, amssymb, amsfonts, yfonts, mathrsfs, float, listings, longtable, mathrsfs, yfonts, tikz-cd, breqn, parskip, enumitem,kantlipsum, mathtools, enumitem, seqsplit, bbm, comment, authblk}
\usepackage[english]{babel}
\usepackage{fancyhdr}
\usepackage{graphicx}
\graphicspath{ {./images/} }
\usepackage[title]{appendix}
\usepackage{csquotes}

\usepackage[
backend=biber,
style=alphabetic,
sorting=nyt,
maxalphanames=10
]{biblatex}
\allowdisplaybreaks

\newtheorem{theorem}{Theorem}[section]
\newtheorem{corollary}[theorem]{Corollary}
\newtheorem{lemma}[theorem]{Lemma}
\newtheorem{remark}[theorem]{Remark}

\title{An Invariance Principle for the 1D KPZ Equation with Neumann Boundary Conditions}

\author[1]{Yiming Tang}
\affil[1]{Department of Mathematics, University of Toronto}

\date{June 2025}

\begin{document}

\newcommand{\eqdis}{\stackrel{d}{=}}
\newcommand{\LL}{\mathcal{L}}
\newcommand{\rr}{\mathbb{R}}
\newcommand{\nn}{\mathbb{N}}
\newcommand{\zz}{\mathbb{Z}}
\newcommand{\YT}[1]{\textcolor{red}{{\sf (Yiming:} {\sl{#1})}}}
\newcommand{\fix}[1]{\textcolor{red}{{\sf (Fix:} {\sl{#1})}}}

\maketitle

\begin{abstract}
    We consider a discrete one-dimensional random interface on the half-space whose height at any positive point is composed of a function of the heights at its two closest neighbours and an independent random noise background. In \cite{adhikari2022invariance}, Adhikari and Chatterjee proved for the full-space model that the height function of such a Markov process converges to the Cole-Hopf solution of the 1D KPZ equation under a parabolic rescaling as the variance of the noise variables goes to zero in the intermediate disorder regime, assuming the dependency of neighboring heights is equivariant, symmetric, and at least six times differentiable in a neighborhood of zero. In this paper, we obtained the same convergence result for the half-space model with a Neumann boundary condition.
\end{abstract}

\section{Introduction}

The 1-dimensional KPZ equation
\begin{equation}
    \partial_t h(x,t)= (\partial_x h(x,t))^2+\partial_x^2 h(x,t) + \xi,
\end{equation}
is a two-dimensional stochastic partial differential equation introduced by Kardar, Parisi, and Zhang \cite{KPZ86} to model a stochastic interface growth process, with $\xi$ being the space-time white noise. In relations with other SPDEs, its partial derivative $\partial_x h(x,t)$ solves the stochastic Burgers equation and $\exp(h(x,t))$ is a mild solution of the stochastic heat equation (SHE) through the Cole-Hopf transformation. The latter was made precise through regularity structures \cite{Hairer_2014}. 
The KPZ equation, as one of the most important models in the KPZ universality class \cite{wu2023kpzequationdirectedlandscape}, is also a scaling limit of many weakly asymmetric interface growth models; this is known as weak KPZ universality. In particular, we focus on half-space models which involve interactions with a boundary. Previous works, including but not limited to the ones cited, have shown the weak universality result in the half space for directed polymers \cite{Wu20} \cite{Par22} \cite{BC23}, asymmetric exclusion processes \cite{CS18} \cite{Par18}, and stochastic six-vertex models \cite{6V2} under various boundary conditions.

In the recent work \cite{adhikari2022invariance}, Adhikari and Chatterjee proved the weak KPZ universality conjecture for a general type of 1D KPZ equation in the absence of a stationary measure.
It is natural for us to ask if we have the same weak universality result in the half-space. The growing interface could be either rewarded or penalized at the boundary. In the main result Theorem \ref{mainresult}, we will show that the height function converges to the log of the mild solution of the SHE with Robin boundary conditions and multiplicative space-time white noise stated as follow
\begin{align*}
    \partial_{t} \mathcal{Z}(x,t)&=\frac{1}{2} \partial_{x}^{2} \mathcal{Z}(x,t)+\sqrt{2 \mu_{2}} \beta \mathcal{Z}(x,t) \xi, \qquad
    \partial_x \mathcal{Z}(0,t)=A\mathcal{Z}(0,t).
\end{align*}
Here, $A, \nu_2$ and $\beta$ are deterministic constant parameters of the SHE and $\xi$ again is the space-time white noise. We will specify how $\nu_2$ and $\beta$ are determined in the next section as we explicitly define our general type of half-space 1D KPZ equation model.

\section{Main result}
Let us denote the background noise in the non-negative half-space by a set of i.i.d. random variables by $(y(x, t))_{x \in \mathbb{N}_{0}, t \in \mathbb{N}_{0}}$, where $\mathbb{N}_0$ is all natural numbers including 0. We assume $y(x,t)$ has a mean zero and a finite moment generating function, denoted by $m$, in a neighborhood of zero. Let $\psi: \mathbb{R}^{2} \rightarrow \mathbb{R}$ be any function satisfying properties: (a) equivariance under constant shifts, i.e. for any $u, v, c \in \mathbb{R}, \psi(u+c, v+c)=\psi(u, v)+c$; (b) symmetry, i.e. for all $u, v \in \mathbb{R}, \psi(u, v)=\psi(v, u)$; (c) regularity, i.e. $\psi$ is at least six times continuously differentiable in a neighborhood of the origin. 

Taking $\beta:=\partial_{1}^{2} \psi(0,0) \neq 0$ and  $\gamma_N:=1-\frac{A}{\sqrt{N}}$ for some constant A, we now consider a growing random surface $f_{N}: \mathbb{N} \times \mathbb{N} \rightarrow \mathbb{R}$ defined through the recursions
\begin{align} \label{half-line model}
    f_{N}(x, t)&=\psi\left(f_{N}(x-1, t-1), f_{N}(x+1, t-1)\right)+N^{-1 / 4} y(x, t) \quad \text{for } x\geq 1,\\
    f_{N}(0, t)&=\psi\left(f_{N}(1, t-1), f_{N}(1, t-1)\right)+N^{-1 / 4} y(0, t)+\beta^{-1}\log\gamma_N,
\end{align}
where $N\in \mathbb{N}$ which we will send to infinity. We denote the initial condition by $f_{N}(x, 0)=\lambda_N(x)$. Here, we assume that $\Lambda_N(x):=\exp(\beta\lambda_N(x))$ is a function that satisfies $\Lambda_N(\sqrt{N}x)$ converging pointwise to a Hölder $1/2^-$ function, denoted by $\mathfrak{Z}_0(x)$, with a convergence rate of $\sup_{x \in \mathbb{N}_0}\|\Lambda_N(\sqrt{N}x)-\mathfrak{Z}_0(x)\| \lesssim N^{-1/4}$. The notation $A(x,t)\lesssim B(x,t)$ means that there exists a deterministic constant $D$, independent of $x, t$ or $N$, such that $A(x,t)\leq DB(x,t)$ for any $(x,t) \in [0, a\sqrt{N}]
 \times [0, bN]$. Further, we assume that $C^{-1} \leq \Lambda_N(x) \leq C$ uniformly on $\mathbb{N}_0$ for some constant $C$.

Here, $\gamma_N$ is the reflection rate at the boundary $x=0$. The boundary is absorbing if $0<\gamma_N<1$, totally reflecting if $\gamma_N=1$, or rewarding if $\gamma_N>1$. $f_{N}(x, t)$ denotes the height of the surface at location $x$ at time $t$. \eqref{half-line model} means that this height is a function of the heights at location $x-1$ and $x+1$ at time $t-1$, when away from the boundary, plus a random noise decaying with $N$.

\begin{theorem}\label{mainresult}
With $m$ being the moment generating function of the noise variable $y(x,t)$, consider a rescaled and renormalized surface growth process $\widetilde{f}_{N}: N^{-1/2}\mathbb{N}_0 \times N^{-1}\mathbb{N}_0 \rightarrow \mathbb{R}$ defined by 
\begin{equation}\label{def2}
\widetilde{f}_{N}(x, t):=f_{N}(\sqrt{N} x, N t)-(V+\beta^{-1}N\log m(N^{-1 / 4} \beta)+N \psi(0,0)) t.
\end{equation}
for some deterministic constant $V$ to be stated explicitly in Remark \ref{V}.

Then extend $\widetilde{f}_{N}$ to $(\mathbb{R}_{\geq 0})^2$ via a linear interpolation. Assuming the initial condition converges fast enough (to be specified in Section 3), $\exp (\beta \widetilde{f}_{N})$ converges in law to a solution $\mathcal{Z}$ of the SHE on the half space with Robin boundary conditions and multiplicative space-time white noise, with respect to the topology of uniform convergence on compact sets on $C((\mathbb{R}_{\geq 0})^2)$.

\end{theorem}

\begin{remark} \label{V}
    For $x \in \mathbb{Z}$ and $t \in \mathbb{N}_0$, let $p(x, t)$ denote the probability that a simple symmetric random walk on $\mathbb{Z}$ starting from $0$ at time $0$, is at $x$ at time $t$. 
    Define $\Delta(x, t):=p(x+1, t)-p(x-1, t)$ and
    $c:=\frac{1}{24} \partial_{1}^{4} \psi(0,0)+\frac{\beta^{3}}{12}.$ Finally, let $\mu_i$ denote the i-th moment of $y(x,t)$, then $V:=c\Big[\sum_{x \in \mathbb{Z}} \sum_{t=0}^{\infty} \Delta(x, t)^{4}(\mu_{4}-\mu_{2}^{2})+\Big(\sum_{x \in \mathbb{Z}} \sum_{t=0}^{\infty} \Delta(x, t)^{2} \mu_{2}\Big)^{2}\Big]$, which turns out to be the same renormalizing constant as in full space result in \cite{adhikari2022invariance}. We will justify the constant in Section \ref{renorm}.

\end{remark}
\begin{remark}
The idea follows from the same approach as \cite{adhikari2022invariance}, through a comparison between the proposed model and the half-space directed polymer model which was known to converge to the Cole-Hopf solution of the half-space KPZ equation with Neumann boundary conditions from \cite{Wu20}. The major difficulty of our work was to verify if the various bounds in \cite{adhikari2022invariance} remain valid in the presence of a reflection parameter and in a lack of complete independence of the noise background when we later mirror the half-space noise variables to the full space. Another difference is that, in Section \ref{renorm}, we are able to prove precisely that the decorrelation is fast through continuous approximations and then we conclude the same result as in \cite{adhikari2022invariance} in the half-space setting; that is, under the stated assumptions, $f_{N}(x, t)$ under parabolic scaling of space and time converges in law to a limit that differs from the Cole-Hopf solution of the 1D KPZ equation on the half space with Neumann boundary condition by a deterministic linear term in $t$.     
\end{remark}

\begin{remark}
    The three conditions we impose on $\psi$ are the same as the full-space model defined in \cite{adhikari2022invariance}. Let us give a brief reasoning on why such assumptions are natural. For weak KPZ universality, we consider models such that the interface growth is driven by lateral growth, certain smoothing mechanisms, and random forcing. Such a model should only care about the relative heights with its neighbors (not necessarily nearest) instead of any absolute heights, and thus condition (a) is reasonable to assume. Condition (b) is assumed for simplicity of the proof. In condition (c), the first two derivatives are essential as the directed polymer model we wish to compare to must have an inverse temperature $\beta=\partial_{1}^{2} \psi(0,0)$. That we further require $\psi \in C^6$ is purely because the proof technique relies on expressing the recursion via Taylor expansions and controlling the leading terms.
\end{remark}

\subsection{Notations}

By considering $\psi_{0}(u, v):=\psi(u, v)-\psi(0,0)$ and combining with the equivariance property of $\psi$ as in \cite{adhikari2022invariance}, we can safely assume $\psi(0,0)=0$ without loss of generality.

From now on, we will fix $N$ and drop the $N$ subscript for the simplicity of writing. We will use $f: \mathbb{N}_0 \times \mathbb{N}_0 \to \mathbb{R}$ to denote the function $f(x, t):=f_{N}(x, t)-t\beta^{-1} \log m\big(N^{-1 / 4} \beta\big).$
Let us also define a function $\phi: \mathbb{R} \rightarrow \mathbb{R}$ by $\phi(u)=\psi\left(\frac{u}{2},-\frac{u}{2}\right)=\psi(u, 0)-\frac{u}{2}$. Then for any $u, v \in \mathbb{R}$, $\psi(u, v) =\phi(u-v)+\frac{u+v}{2}$, meaning that we are able to rewrite the recursion in \eqref{half-line model} to define $f(x,t)$ recursively with $\phi$, which is $C^{6}$ in a neighborhood of zero.

We will also use big "$O$" and small "$o$" notations in the rest of the paper. We say $A(x,t)=O(N^k)$ if there exists some constant $C$ which does not dependent on $N,x$ or $t$ such that $|A(x,t)| \leq CN^k$ for all $N$ large enough. We will specify what $C$ might depend on in the text, as it varies case by case. Similarly, we say $A(x,t)=o(N^k)$ if $|A(x,t)|/N^k$ converges to 0 as $N \to \infty$ uniformly in the given region of $(x,t)$.

\subsection{Acknowledgments}

The author is grateful to Jeremy Quastel for suggesting the problem and many invaluable discussions on both mathematical contents and writing advice. We would also like to thank the anonymous referees for their helpful comments.

\section{Bounds on the Directed Polymer Partition Functions}\label{directedpolymer}

\subsection{Definition}

We will start this section by giving the definition of the half-space polymer model with appropriate boundary conditions as introduced in \cite{Wu20}. We will also give equivalent definitions \eqref{definition2}, \eqref{6.1} and a discrete chaos expansion identity \eqref{twotermexpansion}, which will be used in the proofs of later sections.

Let $RRW(x, t)$ denote the set of all reflected simple symmetric random walk paths on $\mathbb{N}_0$ that start at time 0 and end at $x$ at time $t$.  For $q \in RRW(x,t)$, define $d(q)$ to be $|\{1 \leq i \leq t: q(i)=0\}|$. Therefore, fixing N, we define the partition function  $Z: \mathbb{N}_0 \times \mathbb{N} \rightarrow \mathbb{R}$ as follows.
\begin{equation} \label{definition1}
    Z_N(x, t):= \sum_{q \in RRW(x, t)} 2^{d(q)-t} \Lambda(q(0)) \prod_{i=1}^{t} \frac{\exp \left(\beta N^{-1 / 4} y(q(i), i)\right)}{m(\beta N^{-1 / 4})}  \gamma^{d(q)},
\end{equation}
with $Z_N(x,0):=\Lambda(x)$. Let us use $f^{\text{poly}}_N$ to denote the log partition function, i.e. $f^{\text{poly}}_N=\frac{1}{\beta} \log(Z_N(x, t))$. As pointed out in \cite{adhikari2022invariance} for the full space, $f^{\text{poly}}_N$ also satisfies the recursion defined in \eqref{half-line model} with $\phi^{\text{poly}}(u):=\frac{1}{\beta}\log(\frac{e^{\beta u/2}+e^{-\beta u/2}}{2})$.

 Let $RW(x,t)$ be the set of all possible paths of a simple symmetric random walk on $\mathbb{Z}$ from time 0 that ends at $(x,t)$. For $x \in \mathbb{N}_0$, any path in $RRW(x,t)$ can be viewed as the absolute value of some paths in $RW(x,t)$. Note that $2^{d(q)}$ paths in $RW(x,t)$ correspond to 1 path in $RRW(x,t)$. So equivalently, we have
\begin{equation} \label{definition2}
    Z(x, t):=\frac{1}{2^{t}} \sum_{q \in RW(x, t)} \Lambda(|q(0)|) \gamma^{d(q)} \left(\prod_{i=1}^{t} \frac{\exp \left(\beta N^{-1 / 4} y(|q|(i), i)\right)}{m\left(\beta N^{-1 / 4}\right)}\right).
\end{equation}

Define $\xi(x,t):=\frac{exp(\beta N^{-1/4} y(x,t))}{m(\beta N^{-1/4})}-1$, which can be considered as a multiplicative noise variable with mean zero. Let $\delta_0(x)$ be 1 if $x=0$ and 0 otherwise. We further define $\tilde{\xi}$ by $1+\tilde{\xi}(x,t)=(1+(\gamma-1)\delta_0(x))(1+\xi(x,t))$, so then
\begin{equation}
    \label{6.1}
    Z(x, t):=\frac{1}{2^{t}} \sum_{q \in R W(x, t)} \Lambda(|q(0)|)\prod_{s=1}^{t}[1+\tilde{\xi}(|q(s)|, s)].
\end{equation}

Note that the above expression is valid for $t=0$ as well, adopting the convention that an empty product equals 1. The result in Lemma 6.1 of \cite{adhikari2022invariance} also applies to $\tilde{\xi}$, which states that for all $N$ large enough (depending only on the $\beta$ and the law of the noise variables) and any $p \geq 1$,
$\mathbb{E}|\tilde{\xi}(z, s)|^{2 p} \leq \tilde{C}(p, \beta)N^{-2 p / 4}$,
where $\tilde{C}(p, \beta)$ is a constant that depends on $p, \beta$ and the law of the noise variables (but not on $N$).

If we define $\Gamma(x,t)=\frac{1}{2}(Z(|x+1),t-1)+Z(|x-1|,t-1))$ on $\mathbb{Z} \times \mathbb{N}_0$, it is well-known that (see Lemma 6.2 in \cite{adhikari2022invariance} for example) for any $(x, t) \in \mathbb{Z} \times \mathbb{N}_0$,
\begin{equation} \label{twotermexpansion}
    Z(x, t)=\sum_{z \in \mathbb{Z}}\Lambda(|z|)p(x-z,t)+\sum_{z \in \mathbb{Z}} \sum_{s=1}^{t} p(x-z, t-s) \tilde{\xi}(|z|, s) \Gamma(|z|, s).
\end{equation}

\subsection{Estimates}

In this section, we will state some estimates for $f^{\text {poly }}$ and $Z$, which are almost identical bounds on the full space (Section 4, \cite{adhikari2022invariance}), except that we consider $(x,t) \in [0, a\sqrt{N}] \times [0, bN]$ for Theorem 3.2, Theorem 3.3, and Corollary 3.4. The reason of the discrepancy is due to technical computations on bounding the local time at the boundary by a constant, as demonstrated in detail in Section \ref{sec4}. Even though this rectangular region we consider is strictly smaller than the one considered by \cite{adhikari2022invariance} for the full space, it however corresponds to $[0 , a] \times [0,b]$ after a parabolic rescaling. Thus, we can still obtain the same convergence result on any compact set. 

\begin{lemma} \label{lem4.1}
    Take any fixed $a, b>0$ and $\eta>0$. Suppose in some instance, we have

\begin{equation} \label{4.1}
    \left|\frac{Z(x+1, t)-Z(x-1, t)}{Z(x+1, t)+Z(x-1, t)}\right|<\eta,
\end{equation}
for all $(x, t) \in[1, a N] \times[0, b N]$, then
\begin{equation} \label{4.3}
    \left|f^{\text {poly }}(x+1, t)-f^{\text {poly }}(x-1, t)\right|=\frac{2}{\beta}\left|\frac{Z(x+1, t)-Z(x-1, t)}{Z(x+1, t)+Z(x+1, t)}\right|(1+O(\eta)),
\end{equation}
where the implicit constant above does not depend on $N, x$ or $t$ (but may depend on $a, b$, $\eta$ and $\beta$).

\end{lemma}

\begin{proof}
The proof is identical to the proof of Lemma 4.1 from \cite{adhikari2022invariance}, as the proof only relies on the equality $Z(x,t)=\exp(\beta f^{\text{poly}})$.
\end{proof}

We will postpone proofs of the next three results to Section 5 and 6. Due to the similarity with proofs in \cite{adhikari2022invariance}, we will focus on outlining the differences in our proof.

\begin{theorem} \label{3.2}
    Fix $\epsilon>0, a>0$ and $b>0$. There is an event $\Omega_{L}$ with $\mathbb{P}\left(\Omega_{L}\right)=$ $1-o(1)$ such that on $\Omega_{L}$,
$$
\inf _{0 \leq x \leq a\sqrt{N}, 0 \leq t \leq b N} Z(x, t) \gtrsim N^{-\epsilon} .
$$
\end{theorem}

\begin{theorem} \label{3.3}
    Fix $\epsilon>0, a>0$ and $b>0$. There is an event $\Omega_{U}$ with $\mathbb{P}\left(\Omega_{U}\right)=$ $1-o(1)$, such that on $\Omega_{U}$,
$$
\sup _{1 \leq x \leq a\sqrt{N}, 0 \leq t \leq b N}|Z(x+1, t)-Z(x-1, t)| \lesssim N^{-1 / 4+\epsilon} .
$$
\end{theorem}

As a corollary of these two theorems and Lemma \ref{lem4.1}, we have the following statement.

\begin{corollary}
    Fix $\epsilon>0, a>0$ and $b>0$. Let $\Omega$ be the intersection of the events $\Omega_{L}$ and $\Omega_{U}$ from Theorem \ref{3.2} and Theorem \ref{3.3}. Then $\mathbb{P}(\Omega)=1-o(1)$, and on $\Omega$, we have

$$
\sup _{1 \leq x \leq a\sqrt{N}, 0 \leq t \leq b N}\left|f^{\text {poly }}(x+1, t)-f^{\text {poly }}(x-1, t)\right| \lesssim N^{-1 / 4+\epsilon}
$$
\end{corollary}

\section{Estimates on Exponential Moments of the Local Time} \label{sec4}
In this section we provide  uniform upper and lower bounds on $\mathbb{E}[{Z(x,t)}]$.

\begin{theorem} \label{5.1}
    Fix any $a, b >0$, there exists $C_1, C_2 >0$ such that for all sufficiently large $N$,  $(x,t) \in [0,a\sqrt{N}] \times [0,bN]$, 
    $$C_1 \leq \frac{1}{2^t} \sum_{q \in RW(x,t)} \gamma^{d(q)} \leq C_2.$$
\end{theorem}

\begin{remark}
Note that $\mathbb{E}[{Z(x,t)}]=\frac{1}{2^t}\sum_{q \in RW(x, t)} \Lambda(|q(0)|) \gamma^{d(q)}$. Under our assumption on the initial condition that $C^{-1} \leq \Lambda_N(x) \leq C$ for a  constant $C$ uniform over $N$, it thus suffices to bounds the exponential moment of the local time with a flat initial condition. Combining with Theorem \ref{5.1}, we will obtain  lower and upper bounds uniformly on $[0,a\sqrt{N}] \times [0,bN]$ for $\mathbb{E}[{Z(x,t)}]=\frac{1}{2^t}\sum_{q \in RW(x, t)} \Lambda(|q(0)|) \gamma^{d(q)}$.
\end{remark}

\begin{proof}

Define $TRW(x,t)$ to be the set of simple symmetric random walks on $[0, t]$ starting from $x$, which corresponds to the time reversal of $RW(x,t)$. Let $\tilde{d}(q)$ be the number of visits to 0 of a simple symmetric random walk $q \in TRW(x, t)$ in time $[0, t-1]$. Then the exponential moment of the local time can also be expressed as $\frac{1}{2^{t}}\sum_{q \in TRW(x,t)} \gamma^{\tilde{d}(q)}.$

By Jensen's inequality, as $f(w)=\gamma^w$ is convex, $\frac{1}{2^{t}}\sum_{q \in TRW(x,t)} \gamma^{\tilde{d}(q)} =\mathbb{E}[\gamma^{\tilde{d}}] \geq \gamma^{\mathbb{E}[\tilde{d}]}.$

It is a standard result that $\mathbb{E}[\tilde{d}] \approx \sqrt{t}$, which the $\approx$ sign means there exist non-zero universal constants bounding the ratio of the two sides. So there exist a constant $c>0$ such that $\frac{1}{2^{t}}\sum_{q \in TRW(x,t)} \tilde{d}(q) \leq c\sqrt{t} \leq c\sqrt{bN}$. For $\gamma <1$, the exponential moment of the local time is lower bounded by
$(1-\frac{A}{\sqrt{N}})^{c\sqrt{bN}},$
which converges to $e^{-Ac\sqrt{b}} >0$. For $\gamma\geq 1$, we obtain a trivial lower bound of 1.

To prove the other direction, for any $s \in [1,t]$,
 $|q(s)|-|q(s-1)| = \text{sgn}(q(s-1))(q(s)-q(s-1))+\delta_0(q(s-1))$.

Summing over $s=1$ to $t$, we get
\begin{equation} 
    |q(t)|-|q(0)|=\sum_{s=1}^t \text{sgn}(q(s-1))(q(s)-q(s-1)) + \sum_{s=1}^t \delta_0(q(s-1)).
\end{equation}
Given $|q(0)|=x$, we have
\begin{equation} \label{mean}
    \frac{1}{2^t}\sum_{q \in TRW(x,t)} \gamma^{\tilde{d}(q)}=\frac{1}{2^{t}}\sum_{q \in TRW(x,t)} \gamma^{|q(t)|-x-\sum_{s=1}^t \text{sgn}(q(s-1))(q(s)-q(s-1))}.
\end{equation}

Note that $\gamma^{-x}$ is bounded uniformly over $[0, a\sqrt{N}]$ both above and below by positive constants independent of $N$ as $x=O(\sqrt{N})$ and $\gamma =1-O(N^{-1/2})$. By Cauchy-Schwartz, the rest of (\ref{mean}) is bounded above by 
\begin{equation}
    \Big(\frac{1}{2^{t}}\sum_{q \in TRW(x,t)} \gamma^{2|q(t)|}\Big)^{1/2}\Big(\frac{1}{2^{t}}\sum_{q \in TRW(x,t)} \gamma^{-2\sum_{s=1}^t \text{sgn}(q(s-1))(q(s)-q(s-1))}\Big)^{1/2}.
\end{equation}

For any $x \in [0, a\sqrt{N}]$, and $t \in [0, bN]$, we have
\begin{equation}
        \frac{1}{2^{t}}\sum_{q \in TRW(x,t)} \gamma^{2|q(t)|} =\frac{1}{2^t} \sum_{k=0}^t {{t} \choose {k}} \gamma^{2|x+2k-t|} \leq \frac{1}{2^t} \sum_{k=0}^t {{t} \choose {k}} \left(\gamma^{2(x+2k-t)}+\gamma^{2(-x-2k+t)}\right),
\end{equation}
which equals to $2^{-t}\left(\gamma^{-2}+\gamma^2\right)^t\left(\gamma^{2x}+\gamma^{-2x}\right)$ and therefore it is bounded above uniformly over $(x, t) \in [0, a\sqrt{N}] \times [0, bN]$.

For $\mathbb{E}[\gamma^{-2\sum_{s=1}^t \text{sgn}(q(s-1))(q(s)-q(s-1))}]$, we adopt a standard martingale argument. Take $q \in TRW(x,t)$ and define $\mathcal{F}_s:=\sigma(\{q(i): 1 \leq i \leq s\})$. Then the expectation equals to
\begin{equation}
    \mathbb{E}[\gamma^{-2\sum_{s=1}^{t-1} \text{sgn}(q(s-1))(q(s)-q(s-1))} \mathbb{E}[\gamma^{-2\text{sgn}(q(t-1))(q(t)-q(t-1))}|F_{t-1}]],
\end{equation}
by the tower property.

Note that $\mathbb{E}[\gamma^{-2\text{sgn}(q(t-1))(q(t)-q(t-1))}|F_{t-1}]=\frac{\gamma^2+\gamma^{-2}}{2}(1-\delta_0(q(t-1)))+\delta_0(q(t-1)) \leq \frac{1}{2}\left(\gamma^2+\gamma^{-2}\right)$, which is $1+O(N^{-1})$.

Inductively perform the step above, as $t \leq bN$ and $x \in [0, a\sqrt{N}]$, we would eventually get 
\begin{equation}
    \mathbb{E}\left[\gamma^{-2\sum_{s=1}^t \text{sgn}(q(s-1))(q(s)-q(s-1))}\right] =(1+O(N^{-1}))^t \lesssim 1.
\end{equation}

Combining the inequalities above, we have shown the existence of a constant upper bound for the exponential moment of the local time.

\end{proof}

\section{Proof of Theorem 3.2}

Similar to the whole-space model, we have the following exponential bound.

\begin{theorem}
    Fix $a, b>0$. For any $\theta>0$, there exists some constant $C(\theta)$ that does not depend on $N$ (but may depend on $a$ and $b$) such that for all $(x,t) \in [0, a \sqrt{N}] \times [0, b N]$,
    \begin{equation}
        \mathbb{E}\left[\exp \left(-\theta f^{\text {poly }}(x, t)\right)\right] \leq C(\theta) .
    \end{equation}
\end{theorem}

\begin{proof}
    The proof follows from the idea of Lemma 6.1 from \cite{Chatterjee2021WeakCO}. It remains to show that there are some uniform constants $\alpha$ $C_1$, $C_2$ such that for all $(x,t) \in [0, a\sqrt{N}] \times [0, bN]$ such that $\mathbb{P}(Z(x,t) \geq \alpha) \geq C_1 >0$ and $\mathbb{E}|\Delta Z(x,t)|^2 \geq C_2$.
    
Fix $a,b>0$. Throughout we are taking $(x,t) \in [0,a\sqrt{N}] \times [0, bN]$. Let $\mu:=\frac{m(2\beta N^{-1/4})}{m(\beta N^{-1/4})^2}$. For $q \in RW(x,t)$, define $Z_q(x,t):=\Lambda(|q(0)|)\gamma^{d(q)} \prod_{i=1}^{t} \frac{\exp(\beta N^{-1/4}y(|q|(i),i))}{m(\beta N^{-1/4})} .$ 

 Let $-q$ denote the reflection of $q \in RW(x,t)$ across $0$, which is in $RW(-x,t)$. Define $q \cap r := \{1 \leq i \leq t: q(i)=r(i)\}$. Then for any $q, r \in RW(x,t)$, one has $||q| \cap |r|| \leq |q \cap r| + |(-q) \cap r|$. As $\mu \geq 1$, it follows that $\frac{1}{4^t}\sum_{q, r \in RW(x,t)} \mu^{\left| |q|\cap |r|\right|} \leq \frac{1}{4^t}\sum_{q, r \in RW(x,t)} \mu^{\left| q\cap r\right|} \mu^{\left| (-q)\cap r\right|}$. By Hölder's inequality,

\begin{equation} \label{arg1}
    \begin{aligned}
&\frac{1}{4^t}\sum_{q, r \in RW(x,t)} \mu^{\left| q\cap r\right|} \mu^{\left| (-q)\cap r\right|} \leq \Big(\frac{1}{4^t}\sum_{q, r \in RW(x,t)} \mu^{2\left| q\cap r\right|}\Big)^{1/2}\Big(\frac{1}{4^t}\sum_{q, r \in RW(x,t)} \mu^{2\left| (-q)\cap r\right|}\Big)^{1/2}\\
& \qquad =\mathbb{E}[\mu^{2N_t^{(0)}}]^{1/2} \mathbb{E}[\mu^{2N_t^{(2x)}}]^{1/2},
    \end{aligned}
\end{equation}
where $N_t^{(k)}$ is the number of times that a simple symmetric lazy random walk (each step has 1/4 chance to go to either of its nearest neighbors and 1/2 chance of not moving) hit 0 starting from $k$ in time $(0, t]$. The last equality follows because the difference between two simple symmetric random walks is a simple symmetric lazy random walk. Since $-q$ ends at $-x$ and $r$ ends at $x$, the difference is a lazy walk starting from $2x$.

By Theorem 5.1 in \cite{adhikari2022invariance}, $\mathbb{E}[\mu^{2N_t^{(0)}}] \lesssim 1$. Clearly, $N_{\max(t-x,0)}^{(0)} \leq N_t^{(0)}$ point-wise so 
\begin{equation} \label{arg2}
    \mathbb{E}\left[\mu^{2N_t^{(2x)}}\right] \leq \mathbb{E}\left[\mu^{2N_{\max(t-2x,0)}^{(0)}}\right] \leq \mathbb{E}\left[\mu^{2N_t^{(0)}}\right] \lesssim 1.
\end{equation} Similarly, we can show that for $k=2,4$, $\frac{1}{4^t}\sum_{q, r \in RW(x,t)} \mu^{k\left||q|\cap |r|\right|} \lesssim 1.$

Note that \begin{equation} \label{arg3}
    \mathbb{E}[Z(x,t)^2]  = \mathbb{E}[(\sum_{q \in RW(x,t)} Z_q(x,t))^2] = \frac{1}{4^t}\sum_{q, r \in RW(x,t)} 
|\Lambda(|q(0)|)\Lambda(|r(0)|)|\gamma^{d(q)+d(r)} \mu^{\left| |q|\cap |r|\right|}.
\end{equation}

By Hölder's inequality, it is bounded above by $\Big(\:\frac{1}{4^t}\sum_{q, r \in RW(x,t)} |\Lambda(|q(0)|)\Lambda(|r(0)|)|^2\Big)^{1/2}$ \\ $\big(\:\frac{1}{4^t}\sum_{q, r \in RW(x,t)} \gamma^{4d(q)+4d(r)} \big)^{1/4} \big(\:\frac{1}{4^t}\sum_{q, r \in RW(x,t)}  \mu^{4\left||q|\cap |r|\right|}\big)^{1/4}$. The first term is $O(1)$ by assumption; the second term equals to $(\frac{1}{2^t}\sum_{q \in RW(x,t)} \gamma^{4d(q)} )^{1/2}=O(1)$ by Lemma \ref{5.1}; finally, the last term is $O(1)$ by \eqref{arg2}. Together, there is a constant $C_3$ such that we have $\mathbb{E}[Z(x,t)^2] < C_3$ uniformly. 

Recall Lemma \ref{5.1} gives a constant $C_4>0$ such that $\mathbb{E}[Z(x,t)] \geq C_4$ uniformly. Take $\alpha=\frac{1}{2}C_4$ and by the Paley-Zygmund Inequality, we have
\begin{equation}
    \begin{aligned}
    \mathbb{P}(Z(x,t) \geq \frac{1}{2}C_4) \geq \mathbb{P}(Z(x,t) \geq \frac{1}{2}\mathbb{E}[Z(x,t)]) \geq \frac{1}{4} \frac{\mathbb{E}[Z(x,t)]^2}{\mathbb{E}[Z(x,t)^2]} \geq \frac{C_4^2}{4C_3} > 0.
\end{aligned}
\end{equation}

Next, for any $1 \leq s \leq t$ and $k \in \mathbb{N}_0$, 
\begin{equation}
    \frac{\partial}{\partial y(k,s)} f^{\text{poly}}(x,t)  =\frac{1}{Z(x,t)} \frac{\partial}{\partial y(k,s)} Z(x,t)  =\frac{\beta N^{-1/4}}{2^{t} Z(x,t)} \sum_{q \in RW(x,t), |q(s)|=k} Z_{q}(x,t).
\end{equation}

Let $\nabla f^{\text{poly}}$ denote the gradient of $f^{\text{poly}}$ with respect to $\left(y(k,s)\right)_{1 \leq s \leq t, k \in \mathbb{Z}^{d}}$, then

\begin{align}
|\nabla f^{\text{poly}}(x,t)|^{2} & =\frac{\beta^{2}N^{-1/2}}{4^{ t} Z(x,t)^{2}} \sum_{1 \leq s \leq t} \sum_{k \in \mathbb{N}_0} \; \sum_{\substack{q, r \in RW(x,t), |q(s)|=|r(s)|=k}} Z_q(x,t)Z_r(x,t) \\
& =\frac{\beta^{2}N^{-1/2}}{4^{ t} Z(x,t)^{2}} \sum_{q, r \in RW(x,t)}||q| \cap |r|| Z_q(x,t)Z_r(x,t) .
\end{align}

Hence, $\mathbb{E}|\nabla Z(x,t)|^{2} =\mathbb{E}|Z(x,t) \nabla f^{\text{poly}}(x,t)|^{2}$ can be expanded as 
\begin{equation}
    \frac{1}{4^t}\beta^{2}N^{-1/2} \sum_{q, r \in RW(x,t)}||q| \cap |r|| \mathbb{E}[ Z_q(x,t)Z_r(x,t)],
\end{equation} 

which is bounded by the product of $ \beta^{2}N^{-1/2} ( \frac{1}{4^t}\sum_{q, r \in RW(x,t)}||q| \cap |r||^2)^{1/2}$ and \\ $(\frac{1}{4^t}\sum_{q, r \in RW(x,t)} |\Lambda(q(0))\Lambda(r(0))|^2 \mu^{2||q\cap |r||} \gamma^{2d(q)+2d(r)} )^{1/2}$, using Cauchy-Schwartz. With a similar argument as \eqref{arg3}, we can show that the second term is $O(1)$. It is a standard result (see \cite{adhikari2022invariance} for example) that $\mathbb{E}[(N_t^{(0)})^2] \lesssim N$, so $\mathbb{E}[(N_t^{(2x)})^2] \lesssim N$. Then
\begin{equation}
    \begin{aligned}
        \frac{1}{4^t}\sum_{q, r \in RW(x,t)}||q| \cap |r||^2 \leq \frac{1}{4^t}\sum_{q, r \in RW(x,t)}(|q \cap r|+|(-q) \cap r|)^2 \leq 2\mathbb{E}[(N_t^{(0)})^2] +2\mathbb{E}[(N_t^{(2x)})^2],
    \end{aligned}
\end{equation}
which is also bounded by $N$, up to a constant.

Thus, we conclude that $\mathbb{E}|\nabla Z(x,t)|^{2} \lesssim 1,$ and the rest of the argument in Lemma 6.1 of \cite{Chatterjee2021WeakCO} follows.

\end{proof}

\section{Proof of Theorem 3.3}
For any $x, t \in \mathbb{N}$, let
$\Gamma(x, t):=\frac{1}{2}(Z(x+1, t-1)+Z(x-1, t-1))$ for $x\geq 1$, and $\Gamma(0,t)=Z(1,t-1).$ The following lemma gives a discrete formulation of a pseudo 2-term chaos expansion for the solution of the solution to the half-space stochastic heat equation with Robin boundary conditions. We will use it to compute $Z(x+1)-Z(x-1)$ and then prove Theorem \ref{3.3}. For $0 \leq s \leq t$, $x \geq 1$, define $M(x,s)=\sum_{z \in \mathbb{Z}} \sum_{s=1}^{t} \Delta(x-z, t-s) \tilde{\xi}(|z|, s) \Gamma(|z|, s)$. Then by \eqref{twotermexpansion}, for $t\geq 0$, $x \geq 1$,
$Z(x+1, t)-Z(x-1, t)=\sum_{z \in \mathbb{Z}} \Delta(x-z,t)\Lambda(|z|) +\sum_{z \in \mathbb{Z}} \sum_{s=1}^{t} \Delta(x-z, t-s) \tilde{\xi}(|z|, s) \Gamma(|z|, s),$
and we now denote the second term by $M(x,s)$.

\begin{lemma}
    For any $p \geq 2$, $\left\|M(x,t)\right\|_{L^{p}} \lesssim N^{-1/4}.$
\end{lemma}

\begin{proof}

We split $M(x,t)$ into $M_1(x,t)+M_2(x,t)$, where $M_1$ is deterministic and $M_2$ has mean 0, and bound their norms separately, i.e. $M_1(x,s):=\mathbb{E}\left[M(x,s)\right]$ and $M_2(x,s):=M(x,s)-M_1(x,s)$.

Let $\mathcal{F}_{s}$ be the $\sigma$-algebra generated by all noise variables up to time $s$. Note that for all $s\geq 1$, $\Gamma(z,s) \in F_{s-1}$, so it is independent of $\xi(z,s)$.

We then have $M_1(x,t):= \sum_{r=1}^{t} \sum_{z \in \mathbb{Z}} \Delta(x-z, t-r) \mathbb{E}[\tilde{\xi}(|z|, r)] \mathbb{E}[\Gamma(|z|, r)]$. Recall that $\mathbb{E}[\tilde{\xi}(z, r)]=\gamma-1$ for $z=0$ or $\mathbb{E}[\tilde{\xi}(z, r)]=0$ otherwise. Following from Lemma \ref{5.1},
\begin{equation}
        M_1(x,t):= (\gamma-1) \sum_{r=1}^{t} \Delta(x, t-r)  \mathbb{E}[\Gamma(0, r)] \lesssim -\frac{A}{\sqrt{N}} \sum_{r=1}^{t} \Delta(x, t-r)
        \lesssim N^{-1/2} (1+\sum_{r=1}^{t-1} \Delta(x, r))  
\end{equation}

From Theorem 1.2.1 from \cite{intersections},
\begin{equation}
    \Delta(x, r) = \frac{1}{\sqrt{2\pi r}}\exp\left(-\tfrac{1}{2r}(x+1)^2\right)-\frac{1}{\sqrt{2\pi r}}\exp\left(-\tfrac{1}{2r}(x-1)^2\right)+c_2 O(r^{-2}),
\end{equation}
where $c_2$ is some constant independent of $N$.

As $\sum_{r=1}^{\infty} r^{-1/2}\left(1-\exp\left(2x/r\right)\right)$ converges absolutely. We may conclude that $\|M_1(x,t) \|_{L^p} \lesssim N^{-1/2}.$

Since we have forced $M_2$ to have mean 0, following from the same argument as Section 6 in \cite{adhikari2022invariance} applied to $M_2$, we deduce that $\|M_2(x,t) \|_{L^p} \lesssim N^{-1/4}.$ Then we complete the proof with a triangle inequality.

\end{proof}

By assumptions on the initial conditions in stated in Section 2, we have $\left|\sum_{z \in \mathbb{Z}} \Delta(x-z,t)\Lambda(|z|)\right|=\big|\sum_{z \in \mathbb{Z}}p(x-z,t)\left(\Lambda(|z+1|)-\Lambda(|z-1|)\right)\big|
    = \left|\sum_{z \in \mathbb{Z}}p(x-z,t)\big(\mathfrak{Z}_0(\frac{|z+1|}{\sqrt{N}})-\mathfrak{Z}_0(\frac{|z-1|}{\sqrt{N}})\right.$ \\$\left.+O(N^{-1/4})\big)\right|=O(N^{-1/4}).$

Together, we obtain that for any $p \geq 2$, $\|Z(x+1, t-1)-Z(x-1, t-1)\|_{L^{p}}=O\left(N^{-1 / 4}\right).
$

By Markov's inequality, for any $\epsilon >0$, we have
$\mathbb{P}(|Z(x+1, t-1)-Z(x-1, t-1)| \geq N^{-1/4+\epsilon}) \leq   \frac{\|Z(x+1, t-1)-Z(x-1, t-1)\|^p_{L^{p}}}{(N^{-1/4+\epsilon})^p}
     \lesssim  N^{-\epsilon p}.$
     
Finally, take a union of events $\{ |Z(x+1, t-1)-Z(x-1, t-1)| \geq N^{-1/4+\epsilon}\}$ among $(x,t) \in [1, a\sqrt{N}] \times  [0, bN]$ and define $\Omega_U$ to be its complement, we obtain
\begin{equation}
    \mathbb{P}(\sup_{1 \leq a\sqrt{N}, 0 \leq t\leq bN}|Z(x+1, t-1)-Z(x-1, t-1)| \geq N^{-1/4+\epsilon})
     \lesssim  N^{-\epsilon p+3/2}.
\end{equation}

Taking $p>3/2\epsilon$, we get, with probability $\mathbb{P}(\Omega_U)=1-o(1)$, 
\begin{equation}
    \sup_{1 \leq a\sqrt{N}, 0 \leq t\leq bN}|Z(x+1, t-1)-Z(x-1, t-1)| < N^{-1/4+\epsilon}.
\end{equation}

Combining Lemma 3.1, Theorem 3.2 and 3.3, we have that on the event $\Omega$,
\begin{equation}
    \left(f^{\text {poly }}(x+1, t)-f^{\text {poly }}(x-1, t)\right)^{4}=\frac{16}{\beta^{4}}\left(\frac{Z(x+1, t)-Z(x-1, t)}{Z(x+1, t)+Z(x-1, t)}\right)^{4}+O\left(N^{-1-\epsilon}\right),
\end{equation}
uniformly on a rectangle $[1, a N] \times[0, b N]$.

With the exact same argument in Lemma 7.2, Lemma 7.3, and Corollary 7.4 in \cite{adhikari2022invariance}, we can prove the following two results.
\begin{lemma}
    Fix some $a>0, b>0, \delta>0$ and $\epsilon \in(0,1/100)$. Then there is some event $\Omega_{X}$ such that $\mathbb{P}\left(\Omega_{X}\right)=1-o(1)$ and on $\Omega_{X}$, we have that for all $(x, t) \in[1, a N] \times[0, b N]$,
    \begin{equation}
        \frac{Z(x+1, t)-Z(x-1, t)}{Z(x+1, t)+Z(x-1, t)}=\frac{1}{2} \sum_{z \in \mathbb{Z}} \sum_{t-N^{\epsilon} \leq s \leq t} \tilde{\xi}(|z|, s) \Delta(x-z, t-s)+O\left(N^{-1 / 4-\epsilon / 16}\right),
    \end{equation}
and $K(x,t):=\sum_{z \in \mathbb{Z}} \sum_{t-N^{\epsilon} \leq s \leq t} \tilde{\xi}(|z|, s) \Delta(x-z, t-s) \lesssim N^{-1/4+\delta},$
where the $O$ term is uniform over $(x, t)$ in the above region.
\end{lemma}

\begin{corollary}
    Fix some $a>0, b>0$, and $\epsilon \in(0,1 / 100)$. Then there is some event $\Omega_{4}$ with $\mathbb{P}\left(\Omega_{4}\right)=1-o(1)$ such that on $\Omega_{4}$, we have that for all $(x, t) \in[1, a N] \times[0, b N]$,
\begin{equation}
    \left(\frac{Z(x+1, t)-Z(x-1, t)}{Z(x+1, t)+Z(x-1, t)}\right)^{4}=K(x, t)^{4}+O\left(N^{-1-\epsilon / 32}\right) .
\end{equation} 
\end{corollary}

\section{The Main Argument}
In this section, we will adapt the most important step in proving the result of \cite{adhikari2022invariance} to our half-space model, which estimates the difference between the height function of a directed polymer model and a general height function $f$ as defined in \eqref{half-line model}. Recall the constant $c$ defined in Remark \ref{V} and fix some sufficiently small $\epsilon > 0$, define the "renormalization term"
\begin{equation}\label{y}
    Y(x, t):=\frac{16 c}{\beta^{4}} \sum_{z \in \mathbb{Z}} \sum_{s=1}^{t} p(x-z, t-s) K(|z|, s)^{4} \text {. }
\end{equation}

\begin{theorem} \label{7.1}
    Fix $a, b>0$ and a sufficiently small $\epsilon > 0$. There is an event $\Omega_{2}$ with $\mathbb{P}\left(\Omega_{2}\right)=1-o(1)$, such that on $\Omega_{2}$, we have that for all $(x, t) \in$ $[0, a N] \times[0, b N]$,
$$
f(x, t)=f^{\text {poly }}(x, t)+Y(x, t)+o(1)
$$
where the o(1) term is uniform in $(x, t)$.
\end{theorem}

\begin{proof}
It is easy to check that the same result in Lemma 8.1 in \cite{adhikari2022invariance} follows and we define $\Omega_Y$ to be the event where for all $(x,t) \in [1,a\sqrt{N}] \times [0, bN]$ with $|Y(x+1,t)-Y(x-1,t)| \lesssim N^{-3/2-2\epsilon}$.
Let $\Omega_{2}:=\Omega \cap \Omega_{Y} \cap \Omega_{4}$, where $\Omega$ is from Corollary 3.4 and $\Omega_{4}$ is from Corollary 6.4. Define $\delta(x, t):=f(x, t)-f^{\text {poly }}(x, t)-Y(x, t) .$

The induction argument from Theorem 8.2 in \cite{adhikari2022invariance} follows when we are away from the boundary. So it remains to prove that for all $(x, t)\in D:=\{(x,t) \in \mathbb{N}_0 \times \mathbb{N}_0: x+t \leq \min \{a N, b N\}\}$, having
$|\delta(x, s)| \leq N^{-1-\epsilon / 2} s$ for all $s \leq t-1$ implies $|\delta(0, t-1)| \leq N^{-1-\epsilon / 2} t$ for any $N>N_0$, where $N_0$ is deterministic.

We will work under the assumption that $\Omega_{2}$ holds. Note that $\delta(x, 0)=0$ for all $x$. Fix $(x, t) \in D$. By definition, notice that $Y(0,t)=Y(1,t-1)+\frac{16c}{\beta^4}K(0,t)^4$. Substitute it into the formula of $\delta(x,t)$ to get
\begin{equation}
    |\delta(0,t)|=|\delta(1,t-1)-16c\beta^{-4}K(0,t)^4| \leq |\delta(1,t-1)|+|16c\beta^{-4}K(0,t)^4|.
\end{equation}
From the definition of $K$, as $\Delta(z,s)\xi(|z|,s)$ is an odd function over $z$, $K(0,t)=0$. Hence, we may conclude that 
$|\delta(0,t)|\leq |\delta(1,t-1)| \leq N^{-1-\epsilon/2}t$.
    
\end{proof}

\section{Concentration of the Renormalization Term} \label{renorm}
In this section, we will prove that $Y(x, t)$ asymptotically behaves like a constant multiple of $t$, and will evaluate that constant.

\begin{lemma} \label{lem8.1}
    Fix $a, b>0$ and any sufficiently small $\epsilon>0$. There is an event $\Omega_{Y L L N}$ with $\mathbb{P}\left(\Omega_{Y L L N}\right)=1-o(1)$, such that on $\Omega_{Y L L N}$, we have that for all $(x, t) \in(0, a N] \times[0, b N]$, $|Y(x, t)-\mathbb{E}(Y(x, t))|=o(1)$ where the o(1) term is uniform in $(x, t)$.
\end{lemma} \label{8.1}

     The proof of Lemma \ref{lem8.1} follows exactly the same from the argument as in Lemma 9.1 in \cite{adhikari2022invariance}.

\begin{lemma} \label{8.2}
    Fix $a>0, b>0$, and $\epsilon \in(0,1 / 100)$. Then for any $(x, t) \in \left([0, a ] \times [0, b]\right) \cap \left(N^{-1/2} \mathbb{Z} \times N^{-1} \mathbb{Z}\right)$, we have $\mathbb{E}(Y(\sqrt{N}x, Nt))=V t +o(1)$, where $V$ is the deterministic constant defined in equation (\ref{V}) and the $o(1)$ term is uniformly bounded in $(x, t)$ in the above region.
\end{lemma}
\begin{proof}
Recall that
\begin{equation} \label{Kexpand}
    \begin{aligned}
        K(x, t)&=\frac{1}{2} \sum_{z \in \mathbb{Z}} \sum_{t-N^{\epsilon} \leq s \leq t} \Delta(x-z, t-s) \tilde{\xi}(|z|, s)\\
        &=-\frac{1}{2} \sum_{t-N^{\epsilon} \leq s \leq t} \Delta(x, t-s)
        \frac{A}{\sqrt{N}}+\frac{1}{2} \sum_{z \geq 0} \sum_{t-N^{\epsilon} \leq s \leq t} \hat{\Delta}(x,z, t-s) \xi(z,s),
    \end{aligned}
\end{equation}
where $\hat{\Delta}(x,z, t-s)=\Delta(x-z, t-s)+\Delta(x+z, t-s)$ for $z \neq 0$ and $\hat{\Delta}(x,0, t-s)=\gamma\Delta(x, t-s)$.

We will expand $\mathbb{E}[K(x,t)^4]$ and evaluate each term individually.

Since $m$ is finite in a neighborhood of zero and $\mathbb{E}(y(z, s))=0$, it follows that $m(\theta)=1+O\left(\theta^{2}\right)$ for $\theta$ close to zero. Thus, for any positive integer $k$, as in Lemma 9.2 from \cite{adhikari2022invariance}, we also have
\begin{equation}\label{momentestimate}
    \mathbb{E}[\xi(z, s)^{k}]=\frac{\mathbb{E}\big[\big(e^{\beta N^{-1 / 4} y(z, s)}-m(\beta N^{-1 / 4})\big)^{k}\big]}{m\left(\beta N^{-1 / 4}\right)^{k}}  =\beta^{k} N^{-k / 4} \mu_{k}+o\left(N^{-k / 4}\right).
\end{equation}

This directly tells us the moments of the deterministic term (the first term) in \eqref{Kexpand} is of a smaller magnitude compared to the stochastic term. With the result of Lemma 6.5 in \cite{adhikari2022invariance}, which states $\sum_{z \in \mathbb{Z}} \Delta(z,t)^2 \leq Ct^{-3/2}$, for some constant $C$ independent of $t$. As a result, in the fourth moment of $K(x, t)^{4}$, most terms are bounded by $o\left(N^{-1}\right)$ and we obtained that
\begin{equation}
    \mathbb{E}\left(K(x, t)^{4}\right)=   \mathbb{E}[(\frac{1}{2} \sum_{z \geq 0} \sum_{t-N^{\epsilon} \leq s \leq t} \hat{\Delta}(x,z, t-s) \xi(z,s))^{4}]+o(N^{-1}),
\end{equation}
which can be expanded as $\frac{1}{16} \{\sum_{z \geq 0} \sum_{t-N^{\epsilon} \leq s \leq t} \hat{\Delta}(x, z, t-s)^{4}\left(\mathbb{E}[\xi(z, s)^{4}]-(\mathbb{E}[(\xi(z, s)^{2}])^2\right)$\\ $+\left(\sum_{z \geq 0 } \sum_{t-N^{\epsilon} \leq s \leq t} \hat{\Delta}(x, z, t-s)^{2} \mathbb{E}[\xi(z, s)^{2}]\right)^{2}\}+o(N^{-1})$.

Substituting \eqref{momentestimate} to get $\mathbb{E}\left(K(x, t)^{4}\right)$ can be expressed as
\begin{align} \label{Kestimate}
    \frac{\beta^4}{16N} ((\mu_4-\mu_2^2) \sum_{z \geq 0} \sum_{0 \leq s \leq N^{\epsilon}} \hat{\Delta}(x, z, s)^{4}+\mu_2^2(\sum_{z \geq 0 } \sum_{0 \leq s \leq N^{\epsilon}} \hat{\Delta}(x, z, s)^{2})^{2})+o(N^{-1}).
\end{align}

Fix $x \in (0, a] $, $t \in [0,b]$. To complete the proof, we evaluate $\mathbb{E}(Y(\sqrt{N}x, Nt))$. Substitute the \eqref{Kestimate} to \eqref{y} to get 
\begin{align} 
    &\mathbb{E}(Y(\sqrt{N}x, Nt))=c\sum_{z \in \mathbb{Z}} \sum_{s=1}^{Nt} p(\sqrt{N}x-z,Nt-s)o\left(N^{-1}\right)\\
    & \quad +c\sum_{z \in \mathbb{Z}} \sum_{s=1}^{Nt} p(\sqrt{N}x-z, Nt-s)(\frac{\left(\mu_4-\mu_2^2\right)}{N} \sum_{w \geq 0} \sum_{0 \leq r \leq N^{\epsilon}} \hat{\Delta}(|z|, w, r)^{4}) \label{first_term_y}\\
    & \quad + c\sum_{z \in \mathbb{Z}} \sum_{s=1}^{Nt} p(\sqrt{N}x-z, Nt-s)(\frac{\mu_2^2}{N}\sum_{w \geq 0 } \sum_{0 \leq r \leq N^{\epsilon}} \hat{\Delta}(|z|, w, r)^{2} )^{2}. \label{second_term_y}
\end{align}

Note that the sum of $p(z, s)$ over all $z$ equals 1 for any given $s$ and $t=O(1)$, so the first term is just $o(1)$. We will then deal with the last two non-trivial terms separately. 

As $\hat{\Delta}(z,w,r)=\Delta(z+w,r)+\Delta(z-w,r)$, expanding powers of $\hat{\Delta}$ gives us terms in the form of 
\begin{equation} \label{deltas}
    \sum_{z \in \mathbb{Z}} \sum_{s=1}^{Nt} p(\sqrt{N}x-z, Nt-s) \sum_{w \geq 1} \sum_{0 \leq r \leq N^{\epsilon}} \Delta(|z|+ w, r)^m \Delta(|z|- w, r)^l,
\end{equation} 
for some integers $m,l>0$. The important observation is that $|\Delta(w,s)|$ decreases quickly to 0 when away from $w$. So \eqref{deltas} will become negligible as $\Delta(|z|+ w, r)$ and $\Delta(|z|- w, r)$ move apart. Specifically speaking, we can rewrite \eqref{deltas} as $\sum_{|z|\leq N^{\epsilon}+1} \sum_{s=1}^{Nt} p(\sqrt{N}x-z, Nt-s) \sum_{1\leq |w|\leq N^{\epsilon}+1} \sum_{0 \leq r \leq N^{\epsilon}} \Delta(|z|+ w, r)^m \Delta(|z|- w, r)^l$ since $\Delta(w,r)=0$ whenever $|w|>r+1$. As $\Delta \leq 2$, we finally obtain that \eqref{deltas} is upper bounded by $N^{3\epsilon} \sum_{s=1}^{Nt-1} p(\sqrt{N}x-N^{\epsilon}-1, s)$, up to some constant depending on $m$ and $l$.

Using Theorem 1.2.1 from \cite{intersections} again, approximating the discrete heat kernel by its continuous analog, we get that \eqref{deltas} is bounded by
\begin{align} \label{secondlastbound}
    N^{3\epsilon} \sum_{s=1}^{Nt-1} \Big(\frac{1}{\sqrt{2\pi s}}\exp(-\frac{(\sqrt{N}x-N^{\epsilon}-1)^2}{2s}) + O(s^{-3/2})\Big)  \lesssim  N^{3\epsilon} \sum_{s=1}^{Nt-1} s^{-1/2}
        \lesssim N^{\frac{1}{2}+3\epsilon}.
\end{align}

Applying Lemma 6.5 in \cite{adhikari2022invariance} gives us
$\sum_{r \geq N^{\epsilon}}\sum_{w \in \mathbb{Z}} \Delta( w, r)^4 \leq 4\sum_{r \geq N^{\epsilon}}\sum_{w \in \mathbb{Z}} \Delta( w, r)^2$ $\lesssim  r^{-\epsilon/2} = o(1).$ Combining with \eqref{secondlastbound}, \eqref{first_term_y} becomes
\begin{equation}
N^{-1} \sum_{z \in \mathbb{Z}}\sum_{s=1}^{Nt} p(\sqrt{N}x-z, Nt-s) \Big(\sum_{r \geq 0}\sum_{w \in \mathbb{Z}} \Delta( w, r)^4+o(1)\Big)+o(1)= t \sum_{r \geq 0}\sum_{w \in \mathbb{Z}} \Delta( w, r)^4+o(1).
\end{equation}

Similarly, we can deduce that \eqref{second_term_y} can be expressed as $t (\sum_{r \geq 0}\sum_{w \in \mathbb{Z}} \Delta(w, r)^2)^2 +o(1)$. Altogether, given $V$ in \eqref{V}, we finally obtain that $\mathbb{E}(Y(\sqrt{N}x, Nt))=Vt+o(1).$
\end{proof}

Combining Theorem \ref{7.1}, Lemma \ref{8.1} and Lemma \ref{8.2}, together with the result of which $f^{\text{poly}}$ converges in law to the Cole-Hopf solution of the half-space KPZ equation with Neumann boundary conditions in \cite{Wu20} under parabolic rescaling, we conclude our main result.

%%%%%%%%%%%%%%%%%%%%%%%%%%%%%%%%%%%%%%%%%%%%%%%%%%%%%%%%%%%%%%%%%%%
%%                                                               %%
%% Supplementary Material, if any, should be provided in         %%
%% {supplement} environment  with title and short description.   %%
%%                                                               %%
%%%%%%%%%%%%%%%%%%%%%%%%%%%%%%%%%%%%%%%%%%%%%%%%%%%%%%%%%%%%%%%%%%%

%%%%%%%%%%%%%%%%%%%%%%%%%%%%%%%%%%%%%%%%%%%%%%%%%%%%%%%%%%%%%%%%%%%
%%                                                               %%
%% Use the two commands below for producing your bibliography    %%
%% with bibtex, then comment again the commands and include the  %%
%% content of the .bbl file in this file below the commands.     %%
%%                                                               %%
%%%%%%%%%%%%%%%%%%%%%%%%%%%%%%%%%%%%%%%%%%%%%%%%%%%%%%%%%%%%%%%%%%%

\printbibliography[title={Bibliography}]

\end{document}